\documentclass{article}
\usepackage{amssymb,amsmath,amsthm,graphicx}
\usepackage[all]{xy}
\def\qed{\hfill {\hbox{${\vcenter{\vbox{             
   \hrule height 0.4pt\hbox{\vrule width 0.4pt height 6pt
   \kern5pt\vrule width 0.4pt}\hrule height 0.4pt}}}$}}}

\def\utr{\underline{\triangleright}}
\def\otr{\overline{\triangleright}}

\theoremstyle{definition}
\newtheorem{example}{Example}
\newtheorem{definition}{Definition}
\newtheorem{remark}{Remark}

\date{}

\title{\Large \textbf{Biquandle Image Quivers}}

\author{Geunyoung Kim\footnote{Email: kimg68@mcmaster.ca}
\and
Sam Nelson\footnote{Email: Sam.Nelson@cmc.edu. Partially supported by Simons Foundation collaboration grant 702597}
\and
Puttipong Pongtanapaisan\footnote{Email: 
puttipong$\underline{ }$pongtanpaisan@pitzer.edu}}

\begin{document}
\maketitle

\begin{abstract}
We introduce a new family of weighted quiver-valued invariants of oriented
classical and virtual knots and link, surface-knots and surface-links, and 
immersed surface-knots and surface-links associated to finite biquandles.
We derive several polynomial invariants from these quivers via 
decategorification. We provide examples to illustrate the computation of 
the invariant as well as toy examples showing tables of invariant values 
for small biquandles.
\end{abstract}

\parbox{5.5in} {\textsc{Keywords:} Biquandle homsets, enhancements, 
image enhancement, quivers, surface-links, immersed surface-links,
virtual knots and links

\smallskip

 }

\smallskip

\section{\large\textbf{Introduction}}\label{I}

\textit{Biquandles} are algebraic structures whose axioms encode the
Reidemeister moves of classical and virtual knot theory. Every classical, 
virtual or surface-knot or link $L$ has an associated \textit{fundamental 
biquandle} $\mathcal{B}(L)$ which is a powerful invariant, known to be 
complete up ambient homeomorphism for classical knots. \textit{Quandles}
are a particular type of biquandle.

Given a biquandle $X$, the homset $\mathrm{Hom}(\mathcal{B}(L),X)$ is
a computable invariant whose elements may be represented as colorings of a
diagram of $L$ by elements of $X$ satisfying certain conditions. If $X$ if
finite and $L$ is tame, then the homset is also finite and its cardinality
is an integer-valued invariant called the \textit{biquandle counting 
invariant}, denoted $\Phi_X^{\mathbb{Z}}$. 

An \textit{enhancement} of $\Phi_X^{\mathbb{Z}}$ is a generally stronger
invariant which determines $\Phi_X^{\mathbb{Z}}$. In \cite{CN} the quandle 
homset invariant was categorified to obtain the \textit{quandle coloring 
quiver}, a directed-graph valued invariant which enhances the quandle counting
invariant. In later work such the quiver construction was extended
to categorify other enhancements such as the quandle cocycle invariant
\cite{CN2}, the biquandle bracket invariant \cite{FN}, biquandle module
invariants \cite{IN} and biquandle arrow weight invariants \cite{NS}.

In this paper we extend this construction to one of the earliest 
enhancements, the \textit{image enhancement}, obtaining an infinite
family of categorified invariants of classical knots and links, virtual knots
and links, surface-knots and surface-links, and immersed surface-links.
We decategorify these invariants to obtain new infinite families of 
new polynomial invariants of these objects. The paper is organized as
follows. In Section \ref{BNH} we review the basics of biquandles and their
homset invariants for the various topological objects mentioned. In Section 
\ref{BIQ} we define the biquandle image quiver, categorifying the biquandle
image multiset and polynomial invariants. As an application we define several
new infinite families of invariants and provide a number of toy examples
illustrating the computation of these invariants for some selected small
biquandles. We conclude in Section \ref{Q} with some questions for future 
research.

This paper, including all text, diagrams and computational code,
was produced strictly by the authors without the use of generative AI in any 
form.

\section{\large\textbf{Singular Marked Graph Diagrams}}

We begin with a remark.

\begin{remark}
Throughout this paper we will be considering knotted objects of several types:
oriented classical knots, oriented classical links, oriented virtual knots, 
oriented virtual links, oriented surface-knots, oriented surface-links and 
oriented immersed surface-knots and oriented immersed surface-links. For 
simplicity, we will use the term ``knotted object'' to refer to objects of 
these types.
\end{remark}

Most of the knotted objects of interest in this paper are assumed to
be familiar to the reader. The possible exception is immersed surface-links,
which we will now describe in more detail. 

\begin{definition}
A \textit{singular marked graph diagram} is a planar 4-regular
graph diagram with vertices decorated as \textit{classical crossings}, 
\textit{immersed crossings} and \textit{saddle crossings} as depicted. 

\[\includegraphics{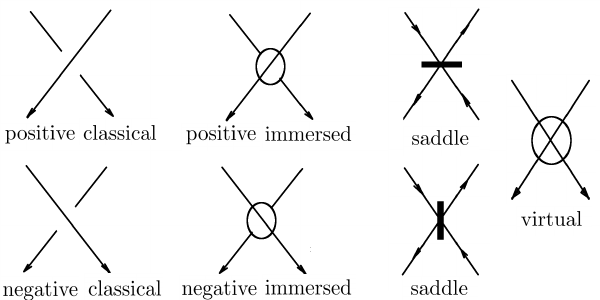}\]

A singular marked graph diagram is \textit{admissible} if the two diagrams 
resulting from resolving all saddle and singular points in both the upper 
smoothing and lower smoothing as shown
\[\includegraphics{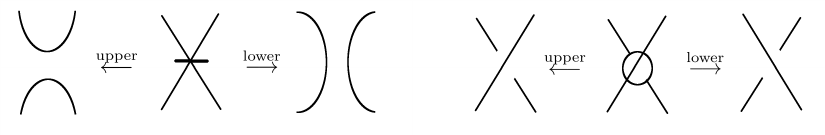}\]
are unlinks. We note that non-admissible diagrams represent cobordisms between
their upper and lower smoothings rather than closed knotted or linked surfaces.
\end{definition}

Two singular marked graph diagrams represent ambient isotopic immersed 
surface-links in $\mathbb{R}^4$ if and only if they are related by the 
\textit{extended Yoshikawa moves} (where the unoriented moves are meant 
to represent all choices of orientation for brevity):
\[\includegraphics{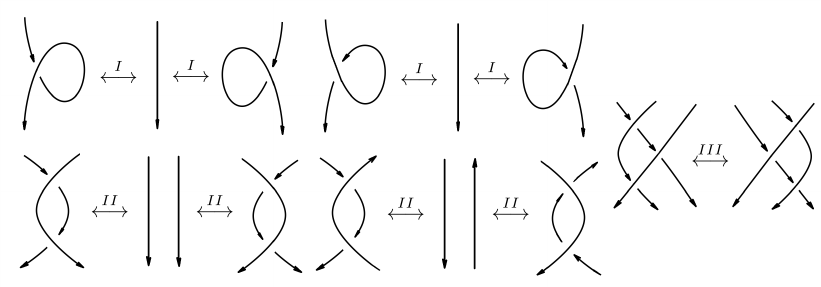}\]
\[\includegraphics{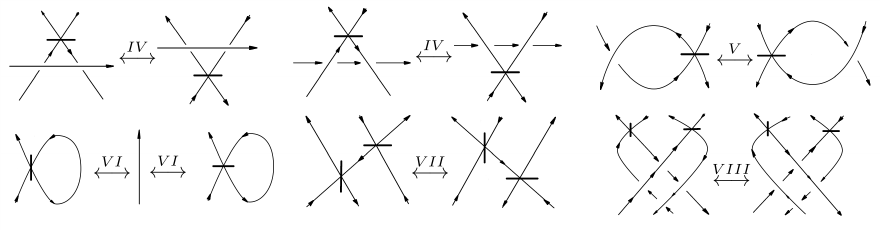}\]
\[\includegraphics{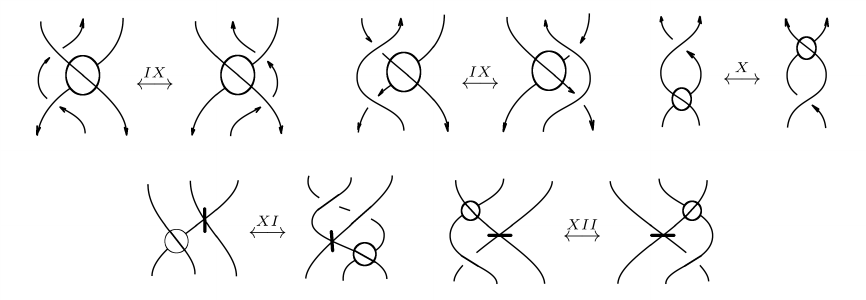}\]
See \cite{HKM,HKM2,J} for more.

\section{\large\textbf{Biquandles and Homsets}}\label{BNH}

In this section we will review the basics of biquandles and their homset
invariants for several types of knotted objects. For more detail see
\cite{EN} and the references therein.

\begin{definition}
A \textit{biquandle} is a set $X$ with two binary operations 
$\utr,\otr:X\times X\to X$ satisfying the axioms
\begin{itemize}
\item[(i)] For all $x\in X$, $x\utr x=x\otr x$,
\item[(ii)] The maps $S:X\times X\to X\times X$ and $\alpha_y,\beta_y:X\to X$
for each $y\in X$ defined by 
\[\alpha_y(x)=x\utr y,\quad \beta_y(x)=x\otr y,\quad \mathrm{and}\quad
S(x,y)=(y\otr x,x\utr y)
\]
are invertible, and
\item[(iii)] The \textit{exchange laws}
\[
\begin{array}{rcl}
(x\utr y)\utr (z\utr y) &=& (x\utr z)\utr (y\otr z) \\
(x\utr y)\otr (z\utr y) &=& (x\otr z)\utr (y\otr z) \\
(x\otr y)\otr (z\otr y) &=& (x\otr z)\otr (y\utr z) 
\end{array}
\]
are satisfied for all $x,y,z\in X$.
\end{itemize}
A biquandle such that $x\otr y=x$ for all $x,y\in X$ is a \textit{quandle}.
\end{definition}

\begin{example}
Standard examples of biquandles include 
\begin{itemize}
\item A group $G$ is a biquandle with $x\utr y=y^{-n}xy^n$ and 
$x\otr y=x$ for any choice of integer $n$,
\item A group $G$ is a biquandle with $x\utr y=y^{-1}xy^{-1}$ and 
$x\otr y=x^{-1}$, and
\item A module $M$ over a commutative unital ring $R$ is 
a biquandle with $x\utr y=tx+(s-t)y$ and $x\otr y=sx$ where $s,t$ are
units in $R$.
\end{itemize}
\end{example}

Biquandles form a category whose objects are biquandles and whose
morphisms are \textit{biquandle homomorphisms}, i.e., maps
$f:X\to Y$ satisfying
\[f(x\utr y)=f(x)\utr f(y)\quad
\mathrm{and}\quad
f(x\otr y)=f(x)\otr f(y)
\]
for all $x,y\in X$. A biquandle homomorphism $f:X\to X$ is an 
\textit{endomorphism}. A subset $S\subset X$ of a biquandle which 
is closed under $\utr$ and $\otr$ is a \textit{sub-biquandle}.

\begin{definition}
Let $L$ be %an oriented knot or link (virtual knot or virtual link, surface-link
%or immersed surface-link) represented by a (knot, virtual knot or marked vertex)
a knotted object represented by a
diagram $D$. Let $X$ be a finite biquandle. Then a \textit{biquandle coloring} 
or an \textit{$X$-coloring} is an assignment of an element of $X$ to each 
semiarc in $D$ such that at each crossing (classical, immersed vertex, 
marked vertex, virtual) we have the following picture:
\[\includegraphics{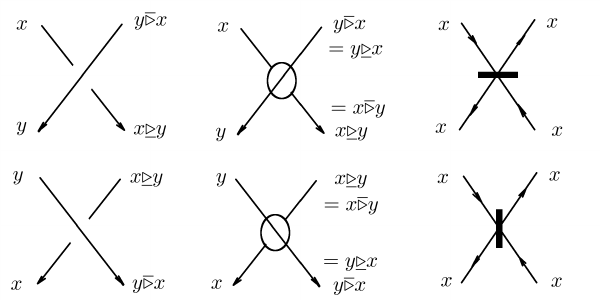}
\raisebox{0.5in}{\includegraphics{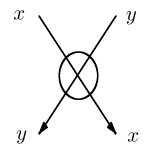}}.\]
We will denote the set of biquandle colorings of $D$ by $X$ as
$\mathrm{Col}(D,X)$.
\end{definition}

\begin{example}
Consider the biquandle $X$ defined by the operation tables
\[
\begin{array}{r|rrrr}
\utr & 1 & 2 & 3 & 4 \\ \hline
1 & 1 & 1 & 2 & 1 \\
2 & 2 & 2 & 1 & 2 \\
3 & 3 & 3 & 3 & 3 \\
4 & 4 & 4 & 4 & 4
\end{array}
\quad
\begin{array}{r|rrrr}
\otr & 1 & 2 & 3 & 4 \\ \hline
1 & 1 & 4 & 2 & 2 \\
2 & 4 & 2 & 1 & 1 \\
3 & 3 & 3 & 3 & 3 \\
4 & 2 & 1 & 4 & 4
\end{array}.
\]
Then the following are examples of $X$-colorings of each of the various
types of knotted objects considered in this paper.
\[
\includegraphics{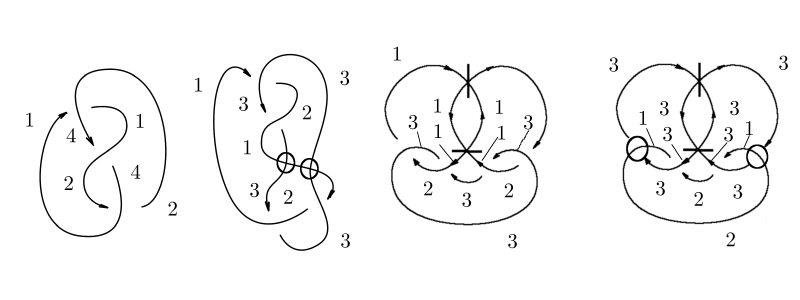}
\]
% 3_1 (1, 1, 2, 2, 4, 4),
% [[-1,2,-3,4],[1,-2,-4,3]] (1, 2, 1, 2, 3, 3, 3, 3)
% [[0, 7, 1, 16, 8], [0, 8, 16, 15, 9], [-1, 1, 12, 2, 13], [-1, 13, 2, 14, 3], [-1, 3, 14, 4, 15], [1, 4, 10, 5, 9], [1, 10, 6, 11, 5], [1, 6, 12, 7, 11]] (4, 4, 1, 1, 1, 4, 4, 4, 4, 2, 2, 1, 2, 2, 4, 4), (4, 4, 2, 2, 2, 4, 4, 4, 4, 1, 1, 2, 1, 1, 4, 4)
% [[0,7,1,16,8],[2,8,16,15,9],[1,1,12,2,13],[-1,13,2,14,3],[-1,3,14,4,15],[1,4,10,5,9],[1,10,6,11,5],[-1,6,12,7,11]],(2, 2, 4, 4, 4, 2, 2, 2, 2, 1, 1, 4, 1, 1, 2, 2)
\end{example}

It is well-known that the set of biquandle colorings of a knotted
object forms an invariant in the case of classical knots and links,
virtual knots and links ans surface-knots and surface-links. The case 
of immersed surface-knots and immersed surface-links was recently shown
in \cite{J}. More precisely, if $D$ and $D'$ are knotted
object diagrams related by Reidemeister/Yoshikawa moves then the sets
of biquandle colorings of $D$ and $D'$ correspond bijectively with 
bijection defined by a sequence of biquandle-colored moves. In particular:
\begin{itemize}
\item The \textit{biquandle counting invariant} of a knotted object
is the cardinality of the set of biquandle colorings; it is a natural 
number-valued invariant.
\item For any invariant $\phi$ of biquandle-colored knots, the multiset
of $\phi$-values over the set of biquandle colorings is an 
\textit{enhancement},
i.e. a new invariant that determines the counting invariant but is in
general a stronger invariant.
\end{itemize}
See \cite{EN} for more.

\begin{definition}
Let $D$ be a diagram representing a knotted object. We obtain a biquandle 
presentation with a generator corresponding to each edge in the underlying
graph of $D$ and a set of relations at each crossing as shown:
\[\scalebox{0.9}{\includegraphics{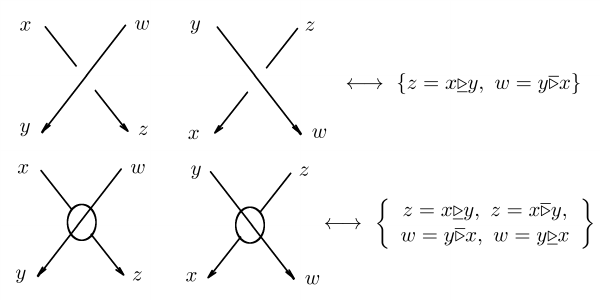}
\includegraphics{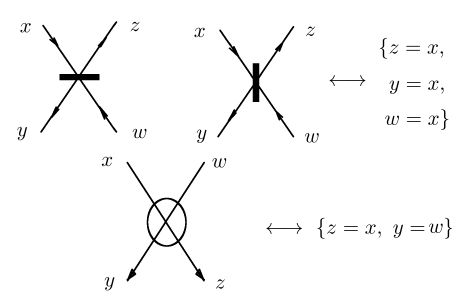}}.\]
We will denote this biquandle as $\mathcal{B}(K)$ where $K$ is the knotted 
object represented by $D$.
\end{definition}

A straightforward check shows that extended Yoshikawa moves on diagrams
induce Tietze moves on the presentation of $\mathcal{B}(K)$; it then
follows that $\mathcal{B}(K)$ is an invariant of the knotted object type in
question. In the cases of classical, virtual and surface-knots and links
$\mathcal{B}(K)$ is called the \textit{fundamental biquandle}. The set of
$X$-colorings of $D$ by a biquandle $X$ can then be identified with the 
homset $\mathrm{Hom}(\mathcal{B}(K),X)$ with the coloring defining an
assignment of image in $X$ to each generator and extending to a biquandle
homomorphism by satisfaction of the coloring rule. In the case
of classical knots, the fundamental biquandle is known to be a complete
invariant up to ambient homeomorphism. See \cite{EN} for more.

\section{\large\textbf{Biquandle Image Quivers}}\label{BIQ}

We begin this section with a definition.

\begin{definition}
Let $X$ be a finite biquandle and let $L$ be a knotted object.
%an oriented classical, virtual, surface- or immersed surface-link.
Let $f:\mathcal{B}(L)\to X$ be a biquandle homomorphism. The \textit{image}
of $f$ is the image of $f$ considered as a map, i.e the sub-biquandle 
\[\mathrm{Im}(f)
=\{x\in X\ |\ x=f(y) \ \mathrm{for \ some }\ y\in\mathcal{B}(L)\}\]
of $X$ generated by the elements which appear in a coloring representing
$f$.
\end{definition}

\begin{example}\label{ex:vk}
Consider the virtual knot $K$ and biquandle $X$ specified by the diagram 
and operation tables 
\[\raisebox{-0.65in}{\includegraphics{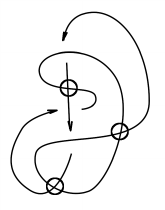}}\quad \quad\quad
\begin{array}{r|rrr}
\utr & 1 & 2 & 3 \\ \hline
1 & 2 & 3 & 1 \\
2 & 3 & 1 & 2 \\
3 & 1 & 2 & 3 \\
\end{array}
\quad
\begin{array}{r|rrr}
\utr & 1 & 2 & 3 \\ \hline
1 & 2 & 2 & 2  \\
2 & 1 & 1 & 1 \\
3 & 3 & 3 & 3
\end{array}.
\]
We compute that there are three $X$-colorings of $K$ by $X$
as shown.
\[\raisebox{-0.5in}{\includegraphics{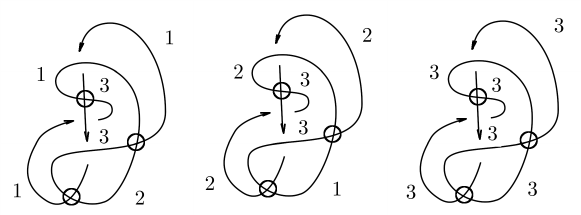}}\]
Then the first two are surjective maps while the third is 
a constant map.
%>>> X[0]
%(((2, 3, 1), (3, 1, 2), (1, 2, 3)), ((2, 2, 2), (1, 1, 1), (3, 3, 3)))
%>>> pdhomlist(gauss2pd(gvknot[(3,1)]),X[0])
%[(1, 3, 3, 1, 1, 2), (2, 3, 3, 2, 2, 1), (3, 3, 3, 3, 3, 3)]
%>>> gvknot[(3,1)]
%[[-1.5, 2, -3.5, -2, 1.5, 3.5]]
\end{example}

We can enhance the biquandle counting invariant by taking the multiset of
cardinalities of image sub-biquandles over the homset.

\begin{definition} (See \cite{EN})
Let $K$ be a knotted object and $X$ a finite biquandle. The 
\textit{image enhancement multiset} of $K$ is the multiset
\[\Phi_X^{\mathrm{Im},M}(K)
=\{|\mathrm{Im}(f)|\ |\ f\in\mathrm{Col}(D,X)\}.\]
\end{definition}

\begin{example}
The virtual knot in Example \ref{ex:vk} has image multiset value
\[\Phi_X^{\mathrm{Im},M}(K)=\{\{1,2,3\},\{1,2,3\},\{3\}\}\]
with respect to the biquandle $X$ in Example \ref{ex:vk}.
\end{example}

In previous work (e.g., \cite{EN}) this multiset invariant was used to define a 
polynomial invariant, namely
\[\Phi_X^{\mathrm{Im},M}(L)=\sum_{f\in\mathrm{Hom}(\mathcal{B}(L),X)} u^{|\mathrm{im}(f)|}.\]

\begin{example}
The virtual knot in Example \ref{ex:vk} has polynomial invariant value
\[\Phi_X^{\mathrm{Im},M}(K)=2u^3+u\]
with respect to the biquandle $X$ in Example \ref{ex:vk}.
\end{example}

We will now categorify this invariant using the biquandle quiver construction
originally introduced in \cite{CN}.

\begin{definition}
Let $K$ be a knotted object,
%n oriented classical, virtual, surface- or immersed surface-link, 
$X$ a finite biquandle and $S\subset \mathrm{Hom}(X,X)$ a 
set of endomorphisms of $X$. Then the \textit{biquandle image quiver} of $K$ 
with respect to $X$ and $S$ is the directed graph $\Phi_{X,S}^{\mathrm{Im}}(K)$
with a vertex for each $f\in\mathrm{Hom}(\mathcal{B}(L),X)$ weighted with
$|\mathrm{Im}(f)|$ and a directed edge from a vertex $f$ to a vertex $g$
whenever $g=\sigma f$ for some $\sigma\in S$.
\end{definition}

\begin{example}
Let $X$ be the biquandle given by the operation tables
\[
\begin{array}{r|rrr}
\utr & 1 & 2 & 3 \\ \hline
1 & 2 & 2 & 2 \\
2 & 1 & 1 & 1 \\
3 & 3 & 3 & 3
\end{array}
\quad
\begin{array}{r|rrr}
\otr & 1 & 2 & 3 \\ \hline
1 & 2 & 2 & 1 \\
2 & 1 & 1 & 2 \\
3 & 3 & 3 & 3
\end{array}
\] 
and consider the Hopf link $L2a1$. Then the maps $\{[3,3,3],[2,1,3]\}$
both represent endomorphisms of $X$, and we have biquandle coloring quiver
\[\includegraphics{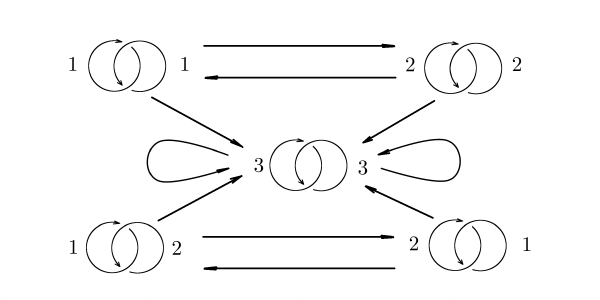}.\]
Then replacing each $X$-coloring in the quiver with a vertex weighted with
the cardinality of its image sub-biquandle, we obtain
\[\Phi_{X,S}^{\mathrm{Im}}(K)=\raisebox{-1in}{\includegraphics{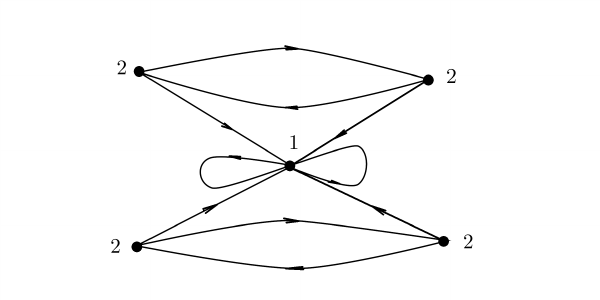}}.\]
\end{example}

Since a quiver is a category, this construction categorifies the image
enhancement multiset and polynomial invariants. The isomorphism class
of the weighted quiver is an invariant which can be used directly to 
distinguish knotted objects. 

As an application, we now define several new 
infinite families of polynomial decategorification invariants.

\begin{definition}
Let $K$ be a knotted object,
%an oriented classical, virtual, surface- or immersed surface-link, 
$X$ a finite biquandle and $S\subset \mathrm{Hom}(X,X)$ a 
set of endomorphisms of $X$. Then the \textit{image in-degree polynomial} 
of $K$ is the polynomial
\[\Phi_{X,S}^{\mathrm{Im},\mathrm{deg}^+}(K)
=\sum_{v\in V} x^{|\mathrm{Im}(v)|}y^{\mathrm{deg}^+(v)}\]
where $V$ is the vertex set of $\Phi_{X,S}^{\mathrm{Im},Q}(K)$.
\end{definition}

\begin{example}
Let $X$ be the biquandle given by the operation tables
\[
\begin{array}{r|rrrrrr}
\utr & 1 & 2 & 3 & 4 & 5 & 6 \\ \hline
1 & 1 & 3 & 2 & 2 & 2 & 2 \\
2 & 3 & 2 & 1 & 3 & 3 & 3 \\
3 & 2 & 1 & 3 & 1 & 1 & 1 \\
4 & 6 & 6 & 6 & 4 & 6 & 5 \\
5 & 5 & 5 & 5 & 6 & 5 & 4\\
6 & 4 & 4 & 4 & 5 & 4 & 6
\end{array}
\quad
\begin{array}{r|rrrrrr}
\otr & 1 & 2 & 3 & 4 & 5 & 6 \\ \hline
1 & 1 & 1 & 1 & 2 & 2 & 2 \\
2 & 2 & 2 & 2 & 3 & 3 & 3 \\
3 & 3 & 3 & 3 & 1 & 1 & 1 \\
4 & 6 & 6 & 6 & 4 & 4 & 4 \\
5 & 5 & 5 & 5 & 5 & 5 & 5 \\
6 & 4 & 4 & 4 & 6 & 6 & 6
\end{array}
\]
and let $S=\{[3,1,2,6,5,4],[2,2,2,2,2,2]\}\subset\mathrm{Hom}(X,X)$.
We then compute via \texttt{python} that the links $L6n1$ and $L7a2$ have
image in-degree polynomial values
\begin{eqnarray*}
\phi_{X,S}^{\mathrm{Im},\mathrm{deg}^+}(L6n1) 
& = & 36x^6y + 18x^4y + xy^{61} + 5xy \\
\phi_{X,S}^{\mathrm{Im},\mathrm{deg}^+}(L6a5) 
& = & 42x^6y + 12x^4y + xy^{61} + 5xy
\end{eqnarray*}
which both specialize to $59y+y^{61}$ at $x=1$, showing that the 
image-enhanced quiver is not determined by the unenhanced quiver.
\end{example}

\begin{example}\label{ex:satoh}
Let $X$ be the conjugation quandle of the symmetric group on four letters.

We then compute via \texttt{python} that the virtual knots $3.7$ and $4.96$ 
have the same counting invariant of $48$. Each strand is labeled with a 
permutation, and we can choose $S\subset \mathrm{Hom}(X,X)$ consisting of a 
single endomorphism that maps an even permutation to the identity 
permutation $(1)$ and maps an odd permutation to the transposition $(12)$. 
Any non-monochromatic diagram of $3.7$ must consist of transpositions as 
labels. On the other hand, a nontrivial coloring of $4.97$ only contains 
3-cycle labels. We see that their in-degree polynomials are distinct:

\begin{eqnarray*}
\phi_{X,S}^{\mathrm{deg}^+}(3.7) 
& = & 30y^5+18y \\
\phi_{X,S}^{\mathrm{deg}^+}(4.96) 
& = & 32y^4+16y
\end{eqnarray*}
Note that by the work of Satoh \cite{SS}, these virtual knot diagrams represent 
ribbon tori in 4-space. Our computation shows that the corresponding 
surfaces are not isotopic.
\end{example}

%\begin{example}
%In \cite{CN}, the authors produce two prime knots $8_{10}$ and $8_{18}$ with the 
%same quandle counting invariant, but non-isomorphic quiver using the dihedral 
%quandle of order 9. One can turn this pair of knots into a pair of immersed 
%surface-links with the same property, which will give a proper enhancement. The 
%readers can use our \texttt{python} code to verify our claim, but we also 
%provide a brief reasoning below.
%
%First, we apply the spinning construction due to Artin. This is formed by 
%removing a small arc $a$ of $K$ to get a knotted arc $K^{\circ}$ in a 3-ball. 
%Taking product with the interval, we obtain a disk in a 4-ball. Doubling such a%disk gives a sphere denoted by spun$(K)$. We can use the result of 
%[CITE Sato and Tanaka] to conclude that the counting invariant of $K$ and 
%spun$(K)$ using dihedral quandles are the same. Furthermore, in any colored 
%diagram of spun$(K)$, there always exist two distinct strands that receive the 
%same color. In fact, these strands come the strand that contained the arc $a$ we 
%removed. Performing a finger move on these two arcs preserve the quandle coloring 
%quiver.
%\end{example}

The classical image enhancement polynomial can be understood as the sum over 
the vertex set in the image quiver of contributions of a variable to the power
of the vertex weight. We can generalize this using the quiver to the get a 
stronger two-variable polynomial invariant.

\begin{definition}
Let $K$ be a knotted object,
%an oriented classical, virtual, surface- or immersed surface-link, 
$X$ a finite biquandle and $S\subset \mathrm{Hom}(X,X)$ a 
set of endomorphisms of $X$. Then the \textit{edge image polynomial} 
of $K$ is the polynomial
\[\Phi_{X,S}^{\mathrm{Im},E}(K)
=\sum_{a\in A} s^{|\mathrm{Im}(S(a))|}t^{|\mathrm{Im}(T(a))|}.\]
where $A$ is the edge set of $\Phi_{X,S}^{\mathrm{Im},Q}(L)$ and $S(a)$, $T(a)$
are the image weights at the source and target of the arrow $a$.
\end{definition}

\begin{example}
Using \texttt{python} code we computed the values of the edge image polynomial
with respect to the biquandle and endomorphism set
\[
\begin{array}{r|rrrr}
\utr & 1 & 2 & 3 & 4 \\ \hline
1 & 1 & 1 & 2 & 1 \\
2 & 2 & 2 & 1 & 2 \\
3 & 3 & 3 & 3 & 3 \\
4 & 4 & 4 & 4 & 4
\end{array}
\quad
\begin{array}{r|rrrr}
\otr & 1 & 2 & 3 & 4 \\ \hline
1 & 1 & 4 & 2 & 2 \\
2 & 4 & 2 & 1 & 1 \\
3 & 3 & 3 & 3 & 3 \\
4 & 2 & 1 & 4 & 4
\end{array}
\quad
S=\{[4,4,3,4],[2,1,3,4],[2,2,2,2]\}
\] 
for the orientable surface-links with ch-index up to 10 from the Yoshikawa 
table \cite{Y}. The results are in the table.
\[
\begin{array}{r|l}
L & \Phi_{X,S}^{\mathrm{Im},E}(L) \\ \hline
2_1 & 12st \\
6^{0, 1}_1 & 4s^4t^4 + 4s^4t^2 + 4s^4t + 4s^2t^2 + 2s^2t + 12st \\
8_1 & 6s^3t^3 + 12s^3t + 12st \\
8^{1, 1}_1 &  4s^4t^4 + 4s^4t^2 + 4s^4t + 4s^2t^2 + 2s^2t + 12st \\
9_1 & 6s^3t^3 + 12s^3t + 12st\\
9^{0, 1}_1 & 4s^4t^4 + 4s^4t^2 + 4s^4t + 4s^2t^2 + 2s^2t + 12st \\
10_1 & 12st \\
10_2 & 6s^3t^3 + 12s^3t + 12st \\
10_3 & 12st \\
10^1_1 & 6s^3t^3 + 12s^3t + 12st \\
10^{0, 1}_1 & 4s^4t^4 + 4s^4t^2 + 4s^4t + 4s^2t^2 + 2s^2t + 12st \\
10^{0,1}_2 & 4s^4t^4 + 4s^4t^2 + 4s^4t + 6s^3t^3 + 12s^3t + 4s^2t^2 + 2s^2t + 12st \\
10^{1, 1}_1 & 4s^4t^4 + 4s^4t^2 + 4s^4t + 4s^2t^2 + 2s^2t + 12st \\
10^{0, 0, 1}_1 & 18s^4t^4 + 18s^4t^2 + 18s^4t + 6s^3t^3 + 12s^3t + 12s^2t^2 + 6s^2t + 12st 
\end{array}.
%[[(2, 1), 12*s*t], [(6, (0, 1), 1), 4*s**4*t**4 + 4*s**4*t**2 + 4*s**4*t + 4*s**2*t**2 + 2*s**2*t + 12*s*t], [(8, 1), 6*s**3*t**3 + 12*s**3*t + 12*s*t], [(8, (1, 1), 1), 4*s**4*t**4 + 4*s**4*t**2 + 4*s**4*t + 4*s**2*t**2 + 2*s**2*t + 12*s*t], [(9, 1), 6*s**3*t**3 + 12*s**3*t + 12*s*t], [(9, (0, 1), 1), 4*s**4*t**4 + 4*s**4*t**2 + 4*s**4*t + 4*s**2*t**2 + 2*s**2*t + 12*s*t], [(10, 1), 12*s*t], [(10, 2), 6*s**3*t**3 + 12*s**3*t + 12*s*t], [(10, 3), 12*s*t], [(10, 1, 1), 6*s**3*t**3 + 12*s**3*t + 12*s*t], [(10, (0, 1), 1), 4*s**4*t**4 + 4*s**4*t**2 + 4*s**4*t + 4*s**2*t**2 + 2*s**2*t + 12*s*t], [(10, (0, 1), 2), 4*s**4*t**4 + 4*s**4*t**2 + 4*s**4*t + 6*s**3*t**3 + 12*s**3*t + 4*s**2*t**2 + 2*s**2*t + 12*s*t], [(10, (1, 1), 1), 4*s**4*t**4 + 4*s**4*t**2 + 4*s**4*t + 4*s**2*t**2 + 2*s**2*t + 12*s*t], [(10, (0, 0, 1), 1), 18*s**4*t**4 + 18*s**4*t**2 + 18*s**4*t + 6*s**3*t**3 + 12*s**3*t + 12*s**2*t**2 + 6*s**2*t + 12*s*t]]
\]

\end{example}

\begin{definition}
Let $K$ be a knotted object,
%an oriented classical, virtual or surface-link, 
$X$ a finite biquandle and $S\subset \mathrm{Hom}(X,X)$ a 
set of endomorphisms of $X$. Then the \textit{image maximal path polynomial}
of $K$ is
\[\Phi_{X,S}^{\mathrm{Im},MP}(K)=\sum_{p\in MP} x^{\sum_{v\in p}w(p)}y^{|p|}\]
where $MP$ is the set of all maximal non-repeating paths in 
the biquandle image quiver, $w(v)$ is the vertex weight at vertex $v$,
and $|p|$ is the length of the path $p$.
\end{definition}

\begin{example}\label{ex:10}
We computed the image maximal path polynomials for the prime classical links 
with up to seven crossings at \cite{KA} with respect to the biquandle given by the
operation tables
\[
\begin{array}{r|rrrrr} 
\utr & 1 & 2 & 3 & 4 & 5 \\ \hline
1 & 2 &3 &1 &2 &2 \\
2 & 3 &1 &2 &1 &1 \\
3 & 1 &2 &3 &3 &3 \\
4 & 4 &4 &4 &5 &5 \\
5 & 5 &5 &5 &4 &4
\end{array}
\quad
\begin{array}{r|rrrrr} 
\utr & 1 & 2 & 3 & 4 & 5 \\ \hline
1 & 2 &2 &2 &2 &2 \\
2 & 1 &1 &1 &1 &1 \\
3 & 3 &3 &3 &3 &3 \\
4 & 5 &5 &5 &5 &5 \\
5 & 4 &4 &4 &4 &4
\end{array}
\]
with set of endomorphisms $S=\{[2,1,3,5,4],[3,3,3,4,5]\}$. The results are in the table:
\[
\begin{array}{r|l}
L & \Phi_{X,S}^{\mathrm{Im},MP}(L) \\ \hline
L2a1 &4x^{12}y^5 + 8x^{10}y^4 \\ 
L4a1 &16x^{30}y^7 + 4x^{12}y^5 + 8x^{10}y^4 \\ 
L5a1 &16x^{30}y^7 + 4x^{12}y^5 + 8x^{10}y^4 \\ 
L6a1 &16x^{30}y^7 + 16x^{12}y^5 + 8x^{10}y^4 \\ 
L6a2 &4x^{12}y^5 + 8x^{10}y^4 \\ 
L6a3 &16x^{12}y^5 + 8x^{10}y^4 \\ 
L6a4 &144x^{30}y^7 + 4x^{12}y^5 + 16x^{10}y^4 \\ 
L6a5 &24x^{30}y^7 + 16x^{12}y^5 + 16x^{10}y^4 \\ 
L6n1 &24x^{30}y^7 + 4x^{12}y^5 + 16x^{10}y^4 \\
\end{array}\quad
\begin{array}{r|l}
L & \Phi_{X,S}^{\mathrm{Im},MP}(L) \\ \hline 
L7a1 &16x^{30}y^7 + 16x^{12}y^5 + 8x^{10}y^4 \\ 
L7a2 &40x^{30}y^7 + 4x^{12}y^5 + 8x^{10}y^4 \\ 
L7a3 &40x^{30}y^7 + 4x^{12}y^5 + 8x^{10}y^4 \\ 
L7a4 &16x^{30}y^7 + 4x^{12}y^5 + 8x^{10}y^4 \\ 
L7a5 &16x^{12}y^5 + 8x^{10}y^4 \\ 
L7a6 &4x^{12}y^5 + 8x^{10}y^4 \\ 
L7a7 &24x^{30}y^7 + 4x^{12}y^5 + 16x^{10}y^4 \\ 
L7n1 &40x^{30}y^7 + 4x^{12}y^5 + 8x^{10}y^4 \\ 
L7n2 &40x^{30}y^7 + 4x^{12}y^5 + 8x^{10}y^4
\end{array}
\]
\end{example}

\begin{example}
The biquandle and endomorphism set in Example \ref{ex:10} also distinguishes
the \textit{Fenn-Rolfsen immersed surface-link} as found in \cite{HKM2,J} 
\[\includegraphics{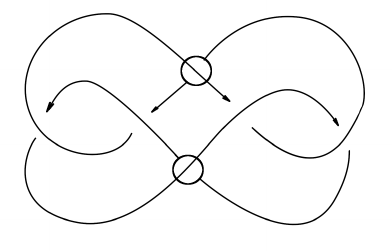}\]
from the unlink of two components with
\[\Phi_{X,S}^{\mathrm{Im},MP}(FR)
=16x^{30}y^7 + 4x^{12}y^5 + 8x^{10}y^4 
\ne 
16x^{30}y^7 + 16x^{12}y^5 + 8x^{10}y^4 =
\Phi_{X,S}^{\mathrm{Im},MP}(U).\]
%>>> FR=[[1,1,8,2,7],[1,12,1,7,6],[-1,9,3,10,4],[-1,4,10,5,11],[2,5,3,6,2],[-2,11,8,12,9]]
%>>> bqcmaxpath(X[7],FR,[[2,1,3,5,4],[3,3,3,4,5]])
%16*x**30*y**7 + 4*x**12*y**5 + 8*x**10*y**4
%>>> bqcmaxpath(X[7],[[1,1,1,2,2],[1,3,3,4,4]],[[2,1,3,5,4],[3,3,3,4,5]])
%16*x**30*y**7 + 16*x**12*y**5 + 8*x**10*y**4

\end{example}

\section{\large\textbf{Questions}}\label{Q}

As with many papers on biquandle- and quiver-based invariants, we caution
the reader that our examples are toy examples, using small biquandles
that are easily accessible via python computation. The true power of this 
infinite family of invariants is expected to lie in the use of larger 
biquandles whose computation currently lie beyond the capacity of our 
python code; therefore, faster computation algorithms are of great interest.

This family of invariants is defined for a wide variety of knotted objects.
It would be of interest to see how the values of the invariant are related 
for similar knotted objects in different categories, expanding on things
like twist-spinning knots or Satoh's tube construction as in Example 
\ref{ex:satoh}.

\bibliography{gk-sn-pp}{}
\bibliographystyle{abbrv}

\bigskip

\noindent
\textsc{Department of Mathematics \& Statistics\\
McMaster University \\
Hamilton, Ontario
\\ Canada}

\smallskip

\noindent
\textsc{Department of Mathematical Sciences \\
Claremont McKenna College \\
850 Columbia Ave. \\
Claremont,  CA 91711}

\smallskip

\noindent
\textsc{Pitzer College \\
1050 Mills Ave \\
Claremont, CA 91711}

\end{document}